\documentclass[11pt]{article}
\usepackage{amsmath,amsthm,amssymb}
\usepackage{subfigure}
\usepackage{graphicx,epstopdf}
\usepackage{euscript,latexsym}
\usepackage[numbers]{natbib}
\bibpunct[, ]{(}{)}{;}{a}{,}{,}
\def\disp{\displaystyle}
\def\dref#1{(\ref{#1})}
\def\crr{\cr\noalign{\vskip2mm}}

\newtheorem{theorem}{Theorem}[section]

\theoremstyle{definition}

\newtheorem{remark}[theorem]{Remark}

\numberwithin{equation}{section}                        

\newcommand{\mylabel}[1]{\label{#1}
            \ifx\undefined\stillediting
            \else \fbox{$#1$}\fi }
\newcommand{\BE}{\begin{equation}}
\newcommand{\EEQ}{\end{equation}}
\newcommand{\rfb}[1]{\mbox{\rm
   (\ref{#1})}\ifx\undefined\stillediting\else:\fbox{$#1$}\fi}

\newfont{\roma}{cmr10 scaled 1200}

\let\oldlabel=\label

\renewcommand{\label}[1]{\leavevmode\smash{\raise 10pt\llap
             {\fbox{\scriptsize#1}}}\oldlabel{#1}}
\newcommand{\mathlabel}[1]{\smash{\raise 9pt\llap
             {\scriptsize(#1)}}\label{#1}}

\renewcommand{\label}[1]{\oldlabel{#1}}
\renewcommand{\mathlabel}[1]{\label{#1}}

\begin{document}

\renewcommand{\thefootnote}{\fnsymbol{footnote}}
\renewcommand{\thefootnote}{\fnsymbol{footnote}}
\newcommand{\footremember}[2]{%
   \footnote{#2}
    \newcounter{#1}
    \setcounter{#1}{\value{footnote}}%
}
\newcommand{\footrecall}[1]{%
    \footnotemark[\value{#1}]%
}
\makeatletter
\def\blfootnote{\gdef\@thefnmark{}\@footnotetext}
\makeatother

\begin{center}
{\Large \bf The Riesz basis property of a class of  \\[0.4ex]
Euler-Bernoulli beam equation}\\[2ex]
Hua-Cheng Zhou
\blfootnote{This work was partially supported by grant no. 800/14
of the Israel Science Foundation.}
\blfootnote{H.-C. Zhou (hczhou@amss.ac.cn) is
with the School of Electrical Engineering, Tel Aviv University,
Ramat Aviv,  69978, Israel.
}
\end{center}

\begin{abstract}
In this paper, we prove that  a sequence of generalized eigenvectors of
a linear unbounded operator associated with an Euler-Bernoulli beam equation
under bending moment boundary feedback forms a Riesz basis for
 the underlying state Hilbert space. As a consequence, the resulting closed-loop
 system is exponentially stable.
\end{abstract}

{\leftskip 10mm {\rightskip 10mm \small \noindent {\bf Key words.}
Riesz basis, exponential stability,
  Euler-Bernoulli beam equation. \bigskip \par}}

\section{Introduction}

Let us consider the following  Euler-Bernoulli beam equation:
\begin{equation}\label{beam-o}
\left\{\begin{array}{l}\disp
w_{tt}(x,t)+w_{xxxx}(x,t)=0,\;0<x<1, \;t\geq0, \crr\disp
w(0,t)=w_x(0,t)=w(1,t)=0, \;t\geq0,\crr\disp
w_{xx}(1,t)=u(t), \;t\geq0,\crr\disp
w(x,0)=w_0(x),\;w_t(x,0)=w_1(x),\;0\leq x\leq1,\crr\disp
y_{m}(t)=w_{xt}(1,t),\;t\geq0.
\end{array}\right.
\end{equation}
In \dref{beam-o}, $w(x,t)$ is the transverse displacement of the beam at time
$t$ and position $x$, $u$ is the input (control) through   bending moment,
$y_m$ is the output signal (measured angular velocity). We are concerned about the following stabilization problem.

{\it Problem.} Given gain $k>0$, does the proportional feedback $u=-ky_m$ make
 the state $(w,w_t)$ of the system  exponentially convergent to zero in the sense that for some $M,\mu>0$,
\begin{equation}\label{beam-Exp}
\int_0^1[|w_t(x,t)|^2+|w_{xx}(x,t)|^2dx\leq Me^{-\mu t}\int_0^1[|w_t(x,0)|^2+|w_{xx}(x,0)|^2dx\;\;\;\;?
\end{equation}
We remark that when boundary condition $w(1,t)=0$ of \dref{beam-o} is replaced by $w_{xxx}(1,t)=0$,
the exponential stability of system under the same proportional feedback $u=-ky_m$ was investigated in \cite{GBZ-YR}.
The motivation studying the above {\it Problem} is to copy with the  Euler-Bernoulli beam
with shear force control matched uncertainties (\cite{ZhouFeng}).
In this paper, we will prove that  the feedback $u=-ky_m$ does make system \dref{beam-o}
exponentially stable in the certain state Hilbert space.

We consider system \dref{beam-o} in the energy Hilbert state space defined by
\begin{equation}
\begin{array}{l}
\mathcal{H}=H_e^2(0,1)\times L^2(0,1), \;
 H_e^2(0,1)=\{\phi\in
H^2(0,1)|\; |\phi(0)=\phi'(0)=0\},
\end{array}
\end{equation}
with the inner product induced norm given by
$$
\|(\phi,\psi)\|^2_{\mathcal{H}}=\int_0^1[|\phi''(x)|^2+|\psi(x)|^2]dx,\; \forall\;
(\phi,\psi)\in\mathcal{H}.
$$
Under this setting and with the feedback $u=-ky_m$, the closed-loop system of \dref{beam-o} can be formulated as
$$
\frac{d}{dt}(w(\cdot,t),w_t(\cdot,t))=A(w(\cdot,t),w_t(\cdot,t)),
$$
where the linear operator is defined as
\begin{equation}\label{A2-def}
\left\{\begin{array}{l}\disp
A(\phi,\psi)=(\psi,-\phi''''),\;\;\forall (\phi,\psi)\in D(A),\;\; \crr\disp
D(A)=\Big\{(\phi,\psi)\in (H^4(0,1)\times H^2_e(0,1))\cap \mathcal{H}: 
                \phi''(1)=-k\psi'(1),\;\phi(1)=0\Big\}.
\end{array}\right.
\end{equation}

\section{Statement and proof of the main results}

Our study in the sequel is focused on the Riesz basis  properties of the
operator $A$. The first main result in this paper can be stated as follows.
\begin{theorem}\label{Lem-A2-exp}
Let $A$ be given by \dref{A2-def}. Then,
there is a sequence of generalized eigenvectors of  $A$  which forms a Riesz basis for
the state space $\mathcal{H}$. Moreover, $A$ generates an exponential stable  $C_0$-semigroup on
$\mathcal{H}$.
\end{theorem}

{\bf Proof.}
The proof is broken into several steps as follows.\\
Step 1. We claim that there is a family of eigenvalues $\{\lambda_n,\overline{\lambda}_n\}$, $\lambda_n=i\tau_n^2$ of  $A$ with the following asymptotic expression:
\begin{equation}\label{taulambda}
\tau_n=(n+{1}/{2})\pi+\mathcal{O}(n^{-1}),\; \;\lambda_n=i(n+1/2)^2\pi^2-\frac{2}{k}+\mathcal{O}(n^{-1}).
\end{equation}
A direct computation shows that
\begin{equation}\label{A2-inv}
\begin{array}{l}\disp
A^{-1}(\phi,\psi)
\hspace{-0.1cm}=\hspace{-0.1cm}\bigg(
\frac{(3x^2-x^3)}{12}\hspace{-0.1cm}\int_0^1\hspace{-0.1cm}(1-\xi)^3\psi(\xi)d\xi
\hspace{-0.1cm}+\hspace{-0.1cm}\frac{(x^3-x^2)}{4}\Big(\hspace{-0.1cm}\int_0^1\hspace{-0.1cm}(1-\xi)\psi(\xi)d\xi\hspace{-0.1cm}-\hspace{-0.1cm}k\phi'(1)\Big)\crr\disp\hspace{3cm}
-\frac{1}{6}\int_0^x(x-\xi)^3\psi(\xi)d\xi,\phi(x)\bigg).
\end{array}
\end{equation}
By the Sobolev embedding theorem, $A^{-1}$ is compact on $\mathcal{H}$, and thus $\sigma(A)$ only consists of
eigenvalues of  $A$.
It is easily seen that $\lambda=i\tau^2\in\sigma(A)$ if and only if there exists $\phi\neq0$ satisfying
\[
\left\{\begin{array}{l}
\phi^{(4)}(x)-\tau^4\phi(x)=0,\crr
\phi(0)=\phi'(0)=\phi(1)=0,\;\;
\phi''(1)=-ik\tau^2\phi'(1)
\end{array}\right.
\]
and the associated egienfunction is $(\phi,\lambda\phi)$. First, the general solution of
\[
\left\{\begin{array}{l}
\phi^{(4)}(x)-\tau^4\phi(x)=0,\crr
\phi(0)=\phi'(0)=0
\end{array}\right.
\]
is of the form
\begin{equation}\label{phi-ex}
\phi(x)=a_1(\cos\tau x-\cosh\tau x)+a_2(\sin\tau x-\sinh\tau x)
\end{equation}
where $a_1$, $a_2$ are constants. Next, by the condition $\phi(1)=0$, we have
$
a_1=(\sin\tau-\sinh\tau),\;\;a_2=-(\cos\tau-\cosh\tau).
$
Substituting this into \dref{phi-ex} gives
\begin{equation}\label{phi-ex1}
\phi(x)=(\sin\tau-\sinh\tau)(\cos\tau x-\cosh\tau x)-(\cos\tau-\cosh\tau)(\sin\tau x-\sinh\tau x).
\end{equation}
The last condition $\phi''(1)=-ik\tau^2\phi'(1)$ yields
\begin{equation}\label{tau-equ}
ik\tau[1-\cos\tau\cosh\tau]+\cosh\tau\sin\tau-\cos\tau\sinh\tau=0,
\end{equation}
which can be re-written asymptotically as
\begin{equation}\label{tau-equ-asy}
\left\{
\begin{array}{l}
\cos\tau=\mathcal{O}(|\tau|^{-1}),\mbox{ or }\crr\disp
\cos\tau=-\frac{1}{ik\tau}[\cos\tau\tanh\tau-\sin\tau]+\mathcal{O}(e^{-{\rm Re}\tau}),
\mbox{ as ${\rm Im}\tau$ bounded ${\rm Re}\tau\to\infty$}.
\end{array}\right.
\end{equation}
By the first equality of \dref{tau-equ-asy}, we get
$
\tau_n=(n+{1}/{2})\pi+\mathcal{O}(n^{-1})
$. Substitute $\tau_n$ into the second equality of \dref{tau-equ-asy} to obtain
$
\mathcal{O}(n^{-1})=-\frac{1}{[ik(n+{1}/{2})\pi]}+\mathcal{O}(n^{-2})
$, and so
\[
\lambda_n=i(n+1/2)^2\pi^2-\frac{2}{k}+\mathcal{O}(n^{-1}).
\]
Step 2.  We claim that there is an eigenfunction $(\phi_n,\lambda_n\phi_n)^\top$ of $A$ corresponding to $\lambda_n=i\tau_n^2$ such that
\[
F_n(x)\hspace{-0.1cm}=\hspace{-0.1cm}\left(\hspace{-0.1cm}\begin{array}{l}
\hspace{-0.1cm}-(-1)^ne^{-(n+1/2)\pi(1-x)}+\cos(n\hspace{-0.1cm}+\hspace{-0.1cm}1/2)\pi x-\sin(n\hspace{-0.1cm}+\hspace{-0.1cm}1/2)\pi x+e^{-(n+1/2)\pi x}\crr
\hspace{-0.1cm}-i(-1)^ne^{-(n+1/2)\pi(1-x)}\hspace{-0.1cm}-\hspace{-0.1cm}i\cos(n\hspace{-0.1cm}+\hspace{-0.1cm}1/2)\pi x\hspace{-0.1cm}+\hspace{-0.1cm}i\sin(n\hspace{-0.1cm}+\hspace{-0.1cm}1/2)\pi x\hspace{-0.1cm}+\hspace{-0.1cm}ie^{-(n+1/2)\pi x}
\end{array}\hspace{-0.2cm}\right)\hspace{-0.1cm}+\mathcal{O}(n^{-1})
\]
and $\lim_{n\to\infty}\|F_n(x)\|_{[L^2(0,1)]^2}=2$,
where $F_n(x)=2\tau_n^{-2}e^{-\tau_n}(\phi_n''(x),\lambda_n\phi_n(x))^\top$. Actually,
let  $(\phi_n,\lambda_n\phi_n)$ be the eigenfunction of $A$ corresponding to $\lambda_n$, where
$\phi_n=\phi(x)$ is defined by \dref{phi-ex} with $\tau=\tau_n$. By \dref{phi-ex1}, we derive
\begin{equation}\label{dotphi}
\tau^{-2}\phi''(x)=(\sin\tau-\sinh\tau)(-\cos\tau x-\cosh\tau x)-(\cos\tau-\cosh\tau)(-\sin\tau x-\sinh\tau x).
\end{equation}
Noticing that by \dref{taulambda}, for any $y>0$ and $0\leq x\leq1$,
$
e^{-\tau_n y}=e^{-(n+1/2)\pi y}+\mathcal{O}(n^{-1}),\;
\sin\tau_n x=\sin(n+1/2)\pi x+\mathcal{O}(n^{-1}),\;
\cos\tau_n x=\cos(n+1/2)\pi x+\mathcal{O}(n^{-1}),
$
and letting $\tau=\tau_n$ in \dref{dotphi}, we obtain
\[
2\tau_n^{-2}\hspace{-0.05cm}e^{-\tau_n}\hspace{-0.04cm}\phi_n''(x)\hspace{-0.1cm}=\hspace{-0.1cm}-(\hspace{-0.05cm}-\hspace{-0.05cm}1)^ne^{-(n\hspace{0cm}+\hspace{0cm}1/2)\pi(1-x)}
\hspace{-0.02cm}+\hspace{-0.02cm}\cos(n+\hspace{-0.05cm}1/2)\pi x-\sin(n+\hspace{-0.05cm}1/2)\pi x+e^{-(n+\hspace{-0.05cm}1/2)\pi x}\hspace{-0.1cm}+\mathcal{O}(n^{-1}\hspace{-0.05cm}).
\]
The estimate for $\phi_n$ is similar, we omit the detail.
By using the Lebesgue's dominated convergence theorem, it is easy to verify $\lim_{n\to\infty}\|F_n(x)\|_{[L^2(0,1)]^2}=2$.\\
Step 3. We claim that the eigenfunctions of $A$  form an Riesz basis for
$\mathcal{H}$. For this purpose, we introduce the following auxiliary operator $A_a$ given by
\begin{equation}\label{A2-aux-def}
\left\{\begin{array}{l}\disp
A_a(\phi,\psi)=(\psi,-\phi''''),\;\;\forall (\phi,\psi)\in D(A_a),\;\; \crr\disp
D(A_a)=\Big\{(\phi,\psi)\in (H^4(0,1)\times H^2_e(0,1))\cap \mathcal{H}:
                \phi''(1)=0,\;\phi(1)=0\Big\}.
\end{array}\right.
\end{equation}
By letting   $k=0$ in \dref{A2-inv}, we know that $A_a$ has compact resolvent.
It is easily to verify that the operator $A_a$ is skew-adjoint in the state space $\mathcal{H}$, i.e.,
$A_a^*=-A_a$ and all eigenvalues of $A_a$ are located on the imaginary axis and there
is a sequence of generalized eigenfunctions
of $A_a$ forming a Riesz basis for $\mathcal{H}$.
Let $\lambda_a=i\omega^2$ be the eigenvalue of $A_a$ and $(\phi_{a},\lambda_{a}\phi_{a})$
be the eigenfunction of $A_a$ corresponding to  $\lambda_a=i\omega^2$.
By letting $k=0$ in \dref{tau-equ}, we obtain
\begin{equation}\label{tau-equ-a}
\cosh\omega\sin\omega-\cos\tau\sinh\omega=0,
\end{equation}
which gives
\[
\omega_n=(n+{1}/{2})\pi+\mathcal{O}(n^{-1}).
\]
Similar to the calculation in Step 2, we obtain that the eigenfunction $(\phi_{an},\lambda_{an}\phi_{an})$
of $A_a$ have the following asymptotical expression:
\[
G_n(x)\hspace{-0.1cm}=\hspace{-0.1cm}\left(\hspace{-0.1cm}\begin{array}{l}
\hspace{-0.1cm}-(\hspace{-0.05cm}-\hspace{-0.05cm}1)^ne^{-(n+1/2)\pi(1-x)}+\cos(n\hspace{-0.1cm}+\hspace{-0.1cm}1/2)\pi x-\sin(n\hspace{-0.1cm}+\hspace{-0.1cm}1/2)\pi x+e^{-(n+1/2)\pi x}\crr
\hspace{-0.1cm}-i(\hspace{-0.05cm}-\hspace{-0.05cm}1)^ne^{-(n+1/2)\pi(1-x)}\hspace{-0.1cm}-\hspace{-0.1cm}i\cos(n\hspace{-0.1cm}+\hspace{-0.1cm}1/2)\pi x\hspace{-0.1cm}+\hspace{-0.1cm}i\sin(n\hspace{-0.1cm}+\hspace{-0.1cm}1/2)\pi x\hspace{-0.1cm}+\hspace{-0.1cm}ie^{-(n+1/2)\pi x}
\end{array}\hspace{-0.2cm}\right)\hspace{-0.1cm}+\mathcal{O}(n^{-1}),
\]
where $G_n(x)=2\omega_n^{-2}e^{-\omega_n}(\phi_{an}''(x),\lambda_n\phi_{an}(x))^\top$.
It is easy to see that $\{(\phi_{an},\lambda_{an}\phi_{an})\}_{n=1}^\infty\cup$ $\{\mbox{conjugates}\}$ is a Riesz basis for $\mathcal{H}$.
It follows  that there is an $N>0$ such that
\begin{equation}\label{quad-close}
\begin{array}{l}\disp
\sum_{n>N}^{\infty}\|F_n-G_n\|_{[L^2(0,1)]^2}=
\sum_{n>N}^{\infty}\|2\tau_n^{-2}e^{-\tau_n}(\phi_n''(x),\lambda_n\phi_n(x))^\top\crr\disp\hspace{0.8cm}
-2\omega_n^{-2}e^{-\omega_n}(\phi_{an}''(x),\lambda_{an}\phi_{an}(x))^\top\|_{[L^2(0,1)]^2}
=\sum_{n>N}^{\infty}\mathcal{O}(n^{-2})<+\infty.
\end{array}
\end{equation}
The same thing is true for conjugates. Therefore,
 operator $A$ has a sequence of eigenfunctions which quadratically closed to
a Riesz basis in the sense of \dref{quad-close}. By \cite[Theorem 1]{GBZ-YR}, we have shown
that the eigenfunctions of $A$  form an Riesz basis for
$\mathcal{H}$.\\
Step 4. We claim that  $A$ generates an exponential stable  $C_0$-semigroup on $\mathcal{H}$.
Since the eigenfunctions of $A$  form an Riesz basis for
$\mathcal{H}$ that is justified by Step 3, the spectrum-determined growth condition holds. In
order to show that $e^{At}$ is a exponential stable semigroup, it suffices to prove that
 ${\rm Re}\lambda<0$   for any $\lambda\in\sigma(A)$. Actually, a simple computation gives
\begin{equation}\label{ReA2-ad}
{\rm Re}\langle A(\phi,\psi),(\phi,\psi)\rangle_{\mathcal{H}}=-k|\psi'(1)|^2\leq0,
\end{equation}
which implies that for any $\lambda\in\sigma(A)$ must
satisfy ${\rm Re}\lambda\leq0$. Since $A^{-1}$   is compact,
we only need to show that there is no eigenvalue on the imaginary
axis. Let $\lambda=i\tau^2\in\sigma(A)$ with
$\tau\in\mathbb{R}^+$ and the corresponding eigenfunction
$(\phi,\psi)^\top\in D(A)$. By \dref{ReA2-ad},
\begin{equation}
{\rm Re}\langle A(\phi,\psi),(\phi,\psi)\rangle_{\mathcal{H}}
={\rm Re}\langle i\tau^2(\phi,\psi),(\phi,\psi)\rangle_{\mathcal{H}}=-k|\psi'(1)|^2=0,
\end{equation}
and hence $\psi'(1)=0$. Furthermore, $A(\phi,\psi)=i\tau^2(\phi,\psi)$ gives that $\psi=i\tau^2\phi$ with $\phi$ satisfying \begin{equation}\label{phi-zerosol}
\left\{\begin{array}{l}
\phi^{(4)}(x)-\tau^4\phi(x)=0,\crr\disp
\phi(0)=\phi'(0)=\phi(1)=\phi'(1)=\phi''(1)=0,
\end{array}\right.
\end{equation}
Now, we show that the above equation admits only zero solution.
For this, we prove that there exists at least one zero of $\phi$ in $(0,1)$.
Actually, by $\phi(0)=\phi(1)=0$, Rolle's theorem yields
$\phi'(\xi_1)=0$ for some $\xi_1\in(0,1)$, which, jointly with $\phi'(0)=\phi'(1)=0$, implies that
$\phi''(\xi_2)=\phi''(\xi_3)=0$ for some $\xi_2\in(0,\xi_1)$,
$\xi_3\in(\xi_1,1)$, and so $\phi'''(\xi_4)=\phi'''(\xi_5)=0$ for some
 $\xi_4\in(\xi_2,\xi_3)$, $\xi_5\in(\xi_3,1)$ by the condition
$\phi''(1)=0$. Thus, there exists a $\xi_6\in(\xi_4,\xi_5)$ such that
$\phi^{(4)}(\xi_6)=0$, which, together with the first equation of
\dref{phi-zerosol}, gives $\phi(\xi_6)=0$.
Next, we prove that if there are $n$ different zeros of $\phi$ in $(0,1)$,
 then there at least $n+1$ number of different zeros of $\phi$ in $(0,1)$.
 Indeed, suppose that $0<\xi_1<\xi_2<\cdots<\xi_n<1$, $\phi(\xi_j)=0,j=1,2,\ldots,n$.
 Since $\phi(0)=\phi(1)=0$, if follows from Rolle's theorem that there exist $\eta_j,j=1,2,\ldots,n+1$, $0<\eta_1<\xi_1<\eta_2<\xi<2\cdots<\xi_n<\eta_{n+1}<1$ such that $\phi'(\eta_j)=0$.
By $\phi'(0)=\phi'(1)=0$, using Rolle's theorem again, there exist $\alpha_j,j=1,2,\ldots,n+2$, $0<\alpha_1<\eta_1<\alpha_2<\eta_2<\cdots<\xi_{n+1}<\alpha_{n+2}<1$ such that $\phi''(\alpha_j)=0$.
It follows from  $\phi''(1)=0$  that there exist $\beta_j,j=1,2,\ldots,n+2$, $\alpha_1<\beta_1<\alpha_2<\beta_2<\cdots<\alpha_{n+2}<\beta_{n+2}<1$ such that $\phi'''(\beta_j)=0$.
 Using Rolle's theorem again, we have $\theta_j,j=1,2,\ldots,n+1$, $\beta_1<\theta_1<\beta_2<\cdots<\beta_{n+1}<\theta_{n+1}<\beta_{n+2}$
 such that $\phi^{(4)}(\theta_j)=0$. Thus, $\phi(\theta_j)=0$, $j=1,2,\ldots,n+1$.
By mathematical induction, there is an infinite number of different zeros $\{x_j\}_{j=1}^\infty$ of $\phi$ in $(0,1)$. Let
$x_0\in[0,1]$ be an accumulation point of $\{x_j\}_{j=1}^\infty$. Obviously, $\phi^{(j)}(x_0)=0$, $j=0,1,2,3$.
Since $\phi$ satisfies the first equation of \dref{phi-zerosol}, by the uniqueness of the solution of linear ordinary different equation,
we have $\phi\equiv0$.
\hfill$\Box$

\begin{remark}
The Riesz basis property of $A$ given by \dref{A2-def} could be used in dealing with disturbance rejection problem considered in \cite{JFF2015auto}
for Euler-Bernoulli beam with shear force control (\cite{ZhouFeng}).
\end{remark}

Now, we state the second result of this paper as follows:
\begin{theorem}
Let system \dref{beam-o} be with the proportional feedback $u=-ky_m$.  
For any initial state $(w(\cdot,0),w_t(\cdot,0))\in\mathcal{H}$,
the resulting closed loop system of \dref{beam-o} admits a unique solution $(w,w_t)\in C(0,\infty;\mathcal{H})$
satisfying \dref{beam-Exp}.
\end{theorem}
{\bf Proof.} The results follows directly from Theorem \ref{Lem-A2-exp}. \hfill$\Box$

%



\end{document}